\documentclass[12 pt]{amsart}
\usepackage{amscd,amsmath,latexsym}
\usepackage[mathscr]{euscript}
\usepackage[all]{xy}
\usepackage{layout}
\usepackage{pb-diagram}
\usepackage{pstricks}

\newtheorem{theorem}{Theorem}[section]
\newtheorem{corollary}[theorem]{Corollary}
\newtheorem{lemma}[theorem]{Lemma}
\newtheorem{proposition}[theorem]{Proposition}
\theoremstyle{definition}
\newtheorem{definition}[theorem]{Definition}
\theoremstyle{definition}

\theoremstyle{remark}

\def\Proof{\medskip\noindent{\bf Proof: }}
 

\def\Z{\mathbb{Z}}

\def\x{\times}
\def\ot{\otimes}

\def\limis{\displaystyle\mathop{\text{lim}^{1}}}
\def\holimi{\displaystyle\mathop{\text{holim}}}
\def\Proof{\medskip\noindent{\bf Proof: }}
\newpsobject{grilla}{psgrid}{subgriddiv=1,griddots=10,gridlabels=6pt}
\def\Proof{\medskip\noindent{\bf Proof: }}

\begin{document}

\title[On the nonexistence of higher twistings]{On the nonexistence of higher twistings}
\author[Jos\'e Manuel G\'omez]{Jos\'e Manuel G\'omez}
\address{Department of Mathematics, University of Michigan, Ann Arbor, MI 48109, USA}
\email{josmago@math.ubc.ca,josmago@umich.edu}
\thanks{The author was supported in part by NSF grant DMS-0244421 and NSF RTG grant 0602191.}

\begin{abstract}
In this note we show that there are no higher twistings for the Borel cohomology theory associated to $G$-equivariant K-theory over a point and for a compact Lie group $G$. Therefore, twistings over a point for this theory are classified by the group $H^{1}(BG,\Z/2)\x H^{3}(BG,\Z)$.
\end{abstract}

\maketitle

\section{Introduction}
The goal of this paper is to show that all the higher twistings for the Borel cohomology theory associated to equivariant K-theory over a point and for a compact Lie group $G$ are trivial.

In general for a non-equivariant multiplicative cohomology theory $E^{*}$, where the multiplication is rigid enough in the sense it is represented by an $E_{\infty}$-ring
spectrum, we can consider local coefficients or twistings. This procedure allows us to construct finer invariants out of the theory $E^{*}$.
The use of local coefficients is a standard tool in algebraic, where for example in the case of singular cohomology they arise naturally in the Serre spectral sequence.

Suppose that $E$ is an $E_{\infty}$-ring
spectrum and let $Z=E_{0}$ be the zero space. $Z$ is an $E_{\infty}$-ring space (see \cite{Mayring} and \cite{MayEinfty} for definitions and  \cite[Corollary 6.6]{MayEinfty}) and  if we write
$Z=\coprod_{\alpha \in \pi_{0}(Z)}Z_{\alpha}$, then $Z_{\otimes}=GL_{1}E=\coprod_{\alpha \in \pi_{0}(Z)^{\times}} Z_{\alpha}$, the space of units, is an infinite loop space by \cite[Corollary 6.8]{MayEinfty}. The space $BZ_{\otimes}$ classifies the twistings for the theory $E^{*}$, this means that for a space $X$ and any map $f:X\to BZ_{\otimes}$, we have a twisting $E^{*}_{f}$ of the theory $E^{*}$ over $X$. The groups $E^{*}_{f}(X)$ and $E^{*}_{f'}(X)$ are isomorphic through a possibly non-canonical isomorphism whenever $f$ and $f'$ are homotopic. In this sense we say that twistings of $E^{*}$ over $X$ are classified by the group $[X,BZ_{\ot}]$. 
As a particular case we can consider non-equivariant K-theory. Twistings for this theory are classified by the spectrum of units $K_{\ot}\simeq \Z/2\times BU_{\otimes}$, where $BU_{\otimes}$ is the space $BU$ with the structure of an H-space corresponding to the tensor product of vector bundles. Thus for a CW-complex $X$, the non-equivariant twistings of complex K-theory over $X$ are classified by the group
\[
[X,BK_{\ot}].
\]
In \cite{Segal}, Segal proved that $BU_{\ot}$ is an infinite loop space and in \cite{Madsen}, Madsen, Snaith and Tornehaveit proved that there is a factorization $BU_{\otimes}=K(\Z,2)\times BSU_{\otimes}$ of the respective spectra. 
We conclude that twistings of K-theory over a space $X$ are classified by homotopy classes of maps  
\[
X\to K(\Z/2,1)\times K(\Z,3)\times BBSU_{\otimes};
\]
that is, we have twistings corresponding to the groups $H^{1}(X,\Z/2)$, 
$H^{3}(X,\Z)$ and $bsu_{\otimes}^{1}(X)=[X,BBSU_{\otimes}]$. We call the twisings corresponding to the group $bsu_{\ot}^{1}(X)$ higher twistings. For the equivariant case the situation is more complicated. We do not know what the ``right'' notion for the spectrum of units for an equivariant spectrum is. In particular, we do not know the group that classifies the more general twistings for equivariant K-theory. However, for the untwisted case, we have the famous Atiyah-Segal completion theorem (see \cite{Atiyahcompletion}). This theorem says that if $G$ is a compact Lie group acting on $X$ a $G$-CW-complex, then we have a natural isomorphism 
\[
K_{G}^{*}(X)\hat{_I}\cong K_{G}^{*}(X\times EG),
\]
where $I$ is the augmentation ideal of the representation ring.

The twistings that we consider here are the twistings of the Borel cohomology theory associated to $G$-equivariant K-theory, where $G$ is a compact Lie group. Hence, for a $G$-CW-complex $X$, the twistings for the Borel cohomology theory associated  to $G$-equivariant K-theory are classified by the group 
\[
H^{1}_{G}(X,\Z/2)\x H^{3}_{G}(X,\Z)\x bsu^{1}_{\otimes}(EG\x_{G} X).
\]
In particular, for the case of a point we obtain
\[
H^{1}(BG,\Z/2)\x H^{3}(BG,\Z)\x bsu^{1}_{\otimes}(BG).
\]
The goal of this paper is to prove that the group $bsu^{1}_{\otimes}(BG)$ vanishes for a compact Lie group, which means that there are no higher twistings for the Borel cohomology theory associated to $G$-equivariant K-theory and over a point.

In what follows we will use the following notation.
\medskip

\noindent{\bf{Notation:}}
We will denote by $k$ the spectrum representing connective complex K-theory and by $K$ the spectrum representing complex K-theory. For a prime $p$ we will denote by $\Z_{p}$ the ring of $p$-adic integers. Given a spectrum $F$ and an abelian group $G$ we can introduce $G$ coefficients on $F$ by considering the spectrum $F_{G}=F\wedge MG$, where $MG$ is a Moore spectrum for the group $G$. Thus in particular we consider $F_{\Z_{p}}$ and $F_{\Z/(p^{k})}$ for a prime number $p$. Also in general for a spectrum $F$ and any integer $n$ we can find the $(n-1)$-connected cover of $F$, which we denote by $F\left\langle n\right\rangle$. This is a spectrum together with a map $F\left\langle n\right\rangle\to F$ that induces an isomorphism $\pi_{k}(F\left\langle n\right\rangle)\stackrel{\cong}{\rightarrow}\pi_{k}F$ for $k\ge n$ and such that $\pi_{k}(F\left\langle n\right\rangle)=0$ for $k<n$. Note that in particular by the periodicity we have $\Sigma^{4}k\cong K\left\langle 4\right\rangle$. 

I would like to thank Professor Robert Bruner for kindly explaining to me that $k^{5}(BG)=0$ for a compact Lie group. This is a crucial result that represents a big part in this work. 

\section{Triviality of $bsu^{1}_{\otimes}(BG)$.} 
In this section we are going to show that the group $bsu^{1}_{\otimes}(BG)$ is trivial. This implies that there are no higher twistings for the Borel cohomology theory associated to equivariant K-theory over a point and for a compact Lie group $G$. In this case, there are only lower twistings; that is, those twistings classified by $H^{1}(BG,\Z/2)\times H^{3}(BG,\Z)$ as in \cite{Atiyah}. 

\begin{definition}
We say that a topological group $G$ satisfies the Atiyah-Segal Completion Theorem if we have that $K^{0}(BG)=R(G)\hat{_I}$ and $K^{1}(BG)=0$, where $I$ is the augmentation ideal of the representation ring $R(G)$.
\end{definition}

Note that by \cite{Atiyahcompletion} it follows that this is true for any compact Lie group.

\begin{lemma}\label{completion}
Let $p$ be a prime number. If $G$ satisfies the Atiyah-Segal Completion Theorem then $K_{\Z_{p}}^{5}(BG)=0$.
\end{lemma}
\Proof
Let us define 
\[
X_{k}=K\wedge M\Z/(p^{k}).
\]
We have maps 
\[
X_{k+1}\to X_{k}
\]
coming from  the maps $\Z/(p^{k+1})\to \Z/(p^{k})$. Consider
\[
X_{\infty}=\holimi_{k\to \infty} X_{k}.
\] 
We will start by showing 
\begin{equation}\label{eq1}
K\wedge M\Z_{p}\cong X_{\infty}.
\end{equation}
We have a map $K\wedge M\Z_{p}\to X_{\infty}$ arising from the canonical maps $\Z_{p}\to \Z/(p^{k})$. Let us show that it induces an isomorphism on homotopy groups. By  \cite[Proposition 6.6]{Adamsbook} there is a short exact sequence 
\begin{equation}\label{eq2}
0\to \pi_{n}(K)\otimes \Z_{p}\to \pi_{n}(K\wedge M\Z_{p})\to \text{Tor}_{1}^{\Z}(\pi_{n-1}(K),\Z_{p})\to 0.
\end{equation}
The group $\text{Tor}_{1}^{\Z}(\pi_{n-1}(K),\Z_{p})$ vanishes, as $\Z_{p}$ is flat as a $\Z$-module. Thus by (\ref{eq2}) we have  
\begin{equation}\label{eq3}
\pi _{n}(K\wedge M\Z_{p})=\left\{ 
\begin{array}{c}
\Z_{p} \text{ if $n$ is even,} \\ 
0 \text{ otherwise.}
\end{array}
\right. 
\end{equation}
On the other hand to compute $\pi_{n}(X_{\infty})$ we have a short exact sequence
\begin{equation}\label{eq4}
0\to \limis_{k\to \infty}\pi_{n+1}(X_{k})\to \pi_{n}(X_{\infty})\to \lim_{k\to \infty}\pi_{n}(X_{k})\to 0.
\end{equation}
Since $\pi_{n}(K)=\Z \text{ or } 0$ according to whether $n$ is even or odd, then by \cite[Proposition 6.6]{Adamsbook} 
$\pi_{n}(X_{k})=\Z/(p^{k}) \text{ or } 0$ depending on the parity of $n$. In any case we have that the map
\[
\pi_{n+1}(X_{k+1})\to \pi_{n+1}(X_{k})
\]
is onto, so the $\text{lim}^{1}$ vanishes. Therefore 
\begin{equation}\label{eq5}
\pi _{n}(X_{\infty})=\left\{ 
\begin{array}{c}
\Z_{p} \text{ if $n$ is even,} \\ 
0 \text{ otherwise}%
\end{array}%
\right. 
\end{equation}%
and the map 
\[
K\wedge M\Z_{p}\to \holimi_{k\to \infty} K\wedge M\Z/(p^{k})=X_{\infty} 
\]
induces isomorphism on $\pi_{*}$. This shows (\ref{eq1}).

Let us show now that 
\begin{equation}\label{eq6}
K_{\Z_{p}}^{5}(BG)=X_{\infty}^{5}(BG)=0.
\end{equation}
To compute this cohomology we consider the short exact sequence
\begin{equation}\label{eq7}
0\to \limis_{k\to \infty}X_{k}^{4}(BG)\to X_{\infty}^{5}(BG)\to \lim_{k\to \infty}X_{k}^{5}(BG)\to 0.
\end{equation}
On the other hand, since $X_{k}=K\wedge M\Z/(p^{k})$, by \cite[Proposition 6.6]{Adamsbook} we have a short exact sequence 
\begin{equation}\label{eq8}
0\to K^{5}(BG)\otimes \Z/(p^{k})\to X_{k}^{5}(BG)\to \text{Tor}^{\Z}_{1}(K^{6}(BG),\Z/(p^{k}))\to 0.
\end{equation}
Since $G$ satisfies the Atiyah-Segal completion theorem, we have that 
\begin{align*}
K^{5}(BG)&=K^{1}(BG)=0\\
K^{6}(BG)&=K^{0}(BG)=R(G)\hat{_I}.
\end{align*}
We know that $R(G)$ is a free, and hence flat $\Z$-module, and $R(G)\hat{_I}$ is a flat $R(G)$-module. By change of basis it follows that $R(G)\hat{_I}$ is a flat $\Z$-module. Therefore from (\ref{eq7}) get that $X_{k}^{5}(BG)=0$.

We also have the exact sequence
\begin{equation}\label{eq9}
0\to K^{4}(BG)\otimes \Z/(p^{k})\to X_{k}^{4}(BG)\to \text{Tor}^{\Z}_{1}(K^{5}(BG),\Z/(p^{k}))\to 0.
\end{equation}
Since $K^{5}(BG)=0$, we conclude from (\ref{eq9}) that $X_{k}^{4}(BG)=K^{4}(BG)\otimes \Z/(p^{k})$. From here we can see that the maps $X_{k+1}^{4}(BG)\to X_{k}^{4}(BG)$ are surjective and thus the $\text{lim}^{1}$ term in the short exact sequence (\ref{eq7}) vanishes. Since the outer terms in that sequence are zero we see that $K_{\Z_{p}}^{5}(BG)=X_{\infty}^{5}(BG)=0$.
\qed

\begin{proposition}\label{Bruner}
Let $G$ be a topological group that satisfies the Atiyah-Segal Completion Theorem. Then $k^{5}(BG)=0$ and $k^{5}_{\Z_{p}}(BG)=0$ for every prime $p$. 
\end{proposition}
\Proof
Both $k^{5}(BG)=0$ and $k^{5}_{\Z_{p}}(BG)=0$ are proved in a similar way with obvious modifications. Thus we will show in detail that $k^{5}_{\Z_{p}}(BG)=0$.

By the previous lemma we have that $K^{5}_{\Z_{p}}(BG)=0$. In general for a spectrum $F$ we have the Atiyah-Hirzebruch spectral sequence.
\[
E_{2}^{r,s}=H^{r}(BG,F^{s}(*))\Longrightarrow F^{r+s}(BG).
\] 
Let us apply this for the cases $F=k_{\Z_{p}}$ and $F=K_{\Z_{p}}$. This way we obtain two spectral sequences $\{E_{n}^{r,s}\}$ and $\{^{1}E_{n}^{ r,s}\}$, respectively.
\begin{align}\label{eq10}
E_{2}^{r,s}=H^{r}(BG,k_{\Z_{p}}^{s}(*))&\Longrightarrow k_{\Z_{p}}^{r+s}(BG),\\
^{1}E_{2}^{r,s}=H^{r}(BG,K_{\Z_{p}}^{s}(*))&\Longrightarrow K_{\Z_{p}}^{r+s}(BG).
\end{align}
For the spectrum $k_{\Z_{p}}$ we know by \cite[Proposition 6.6]{Adamsbook}, that $k_{\Z_{p}}^{n}(*)=\pi_{-n}(k)\otimes \Z_{p}=\Z_{p}$ if $n\le 0$ and even, and $k_{\Z_{p}}^{n}(*)=\pi_{-n}(ku)\otimes \Z_{p}=0$ otherwise. For $K_{\Z_{p}}$ we know that $K^{n}_{\Z_{p}}(*)=\pi_{-n}(K_{\Z_{p}})=\Z_{p}$ if $n$ is even and $\pi_{n}(K_{\Z_{p}})=0$ otherwise. Thus we have that $E_{2}^{r,s}$ is a fourth quadrant spectral sequence with $E_{2}^{r,2s}=H^{r}(BG,\Z_{p})$ for $s\le 0$ and zero otherwise. Similarly, $^{1}E_{2}^{r,2s}=H^{r}(BG,\Z_{p})$ for $s\in \Z$, and zero otherwise. See Figure \ref{1} below.

The spectrum $k$ comes equipped with a map of spectra $k\to K$ inducing an isomorphism on $\pi_{n}$ for $n\ge 0$. By smashing with $MZ_{p}$ we get a map $k_{\Z_{p}}\to K_{\Z_{p}}$ also inducing an isomorphism on $\pi_{n}$ for $n\ge 0$. This map induces a map of spectral sequences $\{E_{n}^{r,s}\}\to \{^{1}E_{n}^{r,s}\}$ as shown in Figure \ref{1}.

We show the result by arguing by contradiction. So assume that $k_{\Z_{p}}^{5}(BG)\ne 0$. We know that $K_{\Z_{p}}^{5}(BG)=0$, and we have a map of spectral sequences $\{E_{n}^{r,s}\}\to \{^{1}E_{n}^{r,s}\}$. Thus the only way that $k_{\Z_{p}}^{5}(BG)\ne 0$ is that one of the differentials that kills elements in total degree 5 in the case $K_{\Z_{p}}$ fails to do so in the case of $k_{\Z_{p}}$. Differentials killing elements in total degree 5 must have source of total degree 4. From Figure \ref{1} we can see at once that the only sources from the $K_{\Z_{p}}$ case of total degree 4 missing in the $k_{\Z_{p}}$ case are $H^{0}(BG,K_{\Z_{p}}^{4}(*))$ and $H^{2}(BG,K_{\Z_{p}}^{2}(*))$. 
We will show that none of these differentials with these sources kill elements of total degree 5 in the case of $K$, from which we deduce that $k_{\Z_{p}}^{5}(BG)=0$.

\begin{center}
\begin{figure}[h]
\centering
\begin{pspicture}(-8,0)(3,6)
\psline{->}(-7.5,0)(-7.5,6)
\psline{->}(-8,3)(-4,3)
\psline{->}(-0.5,0)(-0.5,6)
\psline{->}(-1,3)(3,3)
\psline{->}(-6.5,3)(-5.5,2.5)
\psline{->}(0.5,3)(1.5,2.5)
\pscurve{->}(-3.5,3)(-2.5,3.5)(-1.5,3)
\rput[c](-7.7,6){$s$}
\rput[c](-4,2.8){$r$}
\rput[c](-0.7,6){$s$}
\rput[c](3,2.8){$r$}
\rput[c](-4.5,6){$E_{2}^{r,s}$}
\rput[c](3,6){$\ ^{1}E_{2}^{r,s}$}
\rput[c](-6,2.3){$d^{2}$}
\rput[c](1,2.3){$d^{2}$}
\psdots(-7.5,1)(-7,1)(-6.5,1)(-6,1)(-5.5,1)(-5,1)
\psdots(-7.5,2)(-7,2)(-6.5,2)(-6,2)(-5.5,2)(-5,2)
\psdots(-7.5,3)(-7,3)(-6.5,3)(-6,3)(-5.5,3)(-5,3)
\psdots(-0.5,1)(0,1)(0.5,1)(1,1)(1.5,1)(2,1)
\psdots(-0.5,2)(0,2)(0.5,2)(1,2)(1.5,2)(2,2)
\psdots(-0.5,3)(0,3)(0.5,3)(1,3)(1.5,3)(2,3)
\psdots(-0.5,4)(0,4)(0.5,4)(1,4)(1.5,4)(2,4)
\psdots(-0.5,5)(0,5)(0.5,5)(1,5)(1.5,5)(2,5)
\rput[c](-4,4){$\cdots$}
\rput[c](3,4){$\cdots$}
\end{pspicture}
\caption{Spectral sequences $E, \ ^{1}E$.}
\label{1}
\end{figure}
\end{center}
Let $*$ be the basepoint of $BG$ and consider the sequence of maps  $*\to BG\to *$ factoring the identity $*\to *$. Let us consider now the Atiyah-Hirzebruch spectral sequence applied to the spaces $*$ and $BG$ and the spectrum $K_{\Z_{p}}$. Then we get a spectral sequence $\{^{2}E_{n}^{r,s}\}$
\begin{equation}\label{eq12}
^{2}E_{2}^{r,s}=H^{r}(*,K_{\Z_{p}}^{s}(*))\Longrightarrow K_{\Z_{p}}^{r+s}(*)\\
\end{equation}
and maps $h_{n}^{r,s}:\ ^{1}E_{n}^{r,s}\to\ ^{2}E_{n}^{r,s}$ and $g_{n}^{r,s}:\ ^{2}E_{n}^{r,s}\to\ ^{1}E_{n}^{r,s}$ of spectral sequences such that $h_{n}^{r,s}\circ g_{n}^{r,s}=\text{id}$. (See Figure \ref{Spectral Sequences}.)
\newline
The maps $h_{n}^{r,s}$ and $g_{n}^{r,s}$ and the identity $h_{n}^{r,s}\circ g_{n}^{r,s}=\text{id}$ tell us that all the differentials with source $H^{0}(BG,K_{\Z_{p}}^{4}(*))$ for the spectral sequence $\{^{1}E_{n}^{r,s}\}$ must vanish, as they do for the spectral sequence $\{^{2}E_{n}^{r,s}\}$. 

Now let us study the case of differentials with source $H^{2}(BG,K_{\Z_{p}}^{2}(*))$ for the spectral sequence $\{^{1}E_{n}^{r,s}\}$. We are going to show that all such  differentials are trivial. This is a contradiction and hence the proposition follows.

\begin{center}
\begin{figure}[h]
\centering
\begin{pspicture}(-8,0)(3,6.5)
\psline{->}(-7.5,0)(-7.5,6)
\psline{->}(-8,3)(-4,3)
\psline{->}(-0.5,0)(-0.5,6)
\psline{->}(-1,3)(3,3)
\psline{->}(-7.5,5)(-6.5,4.5)
\psline{->}(-0.5,5)(0.5,4.5)
\pscurve{->}(-3.5,3.2)(-2.5,3.7)(-1.5,3.2)
\pscurve{->}(-1.5,2.8)(-2.5,2.3)(-3.5,2.8)
\rput[c](-2.5,4){$g_{2}^{r,s}$}
\rput[c](-2.5,2){$h_{2}^{r,s}$}
\rput[c](-7.7,6){$s$}
\rput[c](-4,2.8){$r$}
\rput[c](-0.7,6){$s$}
\rput[c](3,2.8){$r$}
\rput[c](-4.5,6){$^{2}E_{2}^{r,s}$}
\rput[c](2.5,6){$^{1}E_{2}^{r,s}$}
\psdots(-7.5,1)
\psdots(-7.5,2)
\psdots(-7.5,3)
\psdots(-7.5,4)
\psdots(-7.5,5)
\psdots(-0.5,1)(0,1)(0.5,1)(1,1)(1.5,1)(2,1)
\psdots(-0.5,2)(0,2)(0.5,2)(1,2)(1.5,2)(2,2)
\psdots(-0.5,3)(0,3)(0.5,3)(1,3)(1.5,3)(2,3)
\psdots(-0.5,4)(0,4)(0.5,4)(1,4)(1.5,4)(2,4)
\psdots(-0.5,5)(0,5)(0.5,5)(1,5)(1.5,5)(2,5)
\rput[c](3,4){$\cdots$}
\end{pspicture}
\caption{Spectral sequences $\ ^{1}E,\ ^{2}E$.}
\label{Spectral Sequences}
\end{figure}
\end{center}

To investigate these differentials we will first study the differentials for the Atiyah-Hirzebruch spectral sequence for the spectrum $K$. So we have a spectral sequence $\{^{3}E_{n}^{r,s}\}$ given by
\begin{equation}\label{eq13}
\ ^{3}E_{2}^{r,s}=H^{r}(BG,K^{s}(*))\Longrightarrow K^{r+s}(BG).
\end{equation}
We are going to show first that all the differentials with source $H^{2}(BG,K^{2}(*))$ vanish. To show this, notice that 
\[
H^{2}(BG,K^{2}(*))=H^{2}(BG,\Z)=[BG,K(\Z,2)],
\] 
and the latter is in a one to one correspondence with isomorphism classes of complex line bundles over $BG$, so every element in $H^{2}(BG,K^{2}(*))$ is the first Chern class of a complex line bundle over $BG$. Let $\alpha\in H^{2}(BG,K^{2}(*))$. Then we can find a map $f:BG\to K(\Z,2)$ such that $\alpha=f^{*}(c_{1}(\gamma_{1}))=c_{1}(f^{*}{\gamma_{1}})$ , where $\gamma_{1}$ is universal line bundle over $K(\Z,2)$. Let ${^{4}E_{n}^{p,q}}$ be the Atiyah-Hirzebruch spectral sequence of the space $K(\Z,2)\simeq \mathbb{C}P^{\infty}$ corresponding to the spectrum $K$, so that
\begin{equation}\label{eq14}
^{4}E_{2}^{r,s}=H^{r}(K(\Z,2),K^{s}(*))\Longrightarrow K^{r+s}(K(\Z,2)).
\end{equation}
The $\ ^{4}E_{2}$-term of this spectral sequence only has terms in the even components and hence the sequence collapse and all the higher differentials are zero.

The map $f$ gives a map of spectral sequences $f_{n}^{r,s}:\ ^{4}E_{n}^{r,s}\to\ ^{3}E_{n}^{r,s}$.
By construction we have that $f_{2}^{2,2}(c_{1}(\gamma_{1}))=\alpha$. Since all the differentials on $\{^{4}E_{n}^{r,s}\}$ are zero it follows that $\alpha$ vanishes on all the differentials. Since $\alpha$ was arbitrary we see that all the differentials on the spectral sequence $\{^{3}E^{p,q}_{n}\}$ with source $H^{2}(BG,K^{2}(*))$ must vanish.

\begin{center}
\begin{figure}[ht]
\centering
\begin{pspicture}(-8,0)(3,6)
\psline{->}(-7.5,0)(-7.5,6)
\psline{->}(-8,3)(-4,3)
\psline{->}(-0.5,0)(-0.5,6)
\psline{->}(-1,3)(3,3)
\psline{->}(-7.5,5)(-6.5,4.5)
\psline{->}(-0.5,5)(0.5,4.5)
\pscurve{->}(-3.5,3.2)(-2.5,3.7)(-1.5,3.2)
\rput[c](-2.5,4){$f_{2}^{r,s}$}
\rput[c](-7.7,6){$s$}
\rput[c](-4,2.8){$r$}
\rput[c](-0.7,6){$s$}
\rput[c](3,2.8){$r$}
\rput[c](-4.5,6){$\ ^{4}E_{2}^{r,s}$}
\rput[c](2.5,6){$\ ^{3}E_{2}^{r,s}$}
\psdots(-7.5,1)(-6.5,1)(-5.5,1)
\psdots(-7.5,2)(-6.5,2)(-5.5,2)
\psdots(-7.5,3)(-6.5,3)(-5.5,3)
\psdots(-7.5,4)(-6.5,4)(-5.5,4)
\psdots(-7.5,5)(-6.5,5)(-5.5,5)
\psdots(-0.5,1)(0,1)(0.5,1)(1,1)(1.5,1)(2,1)
\psdots(-0.5,2)(0,2)(0.5,2)(1,2)(1.5,2)(2,2)
\psdots(-0.5,3)(0,3)(0.5,3)(1,3)(1.5,3)(2,3)
\psdots(-0.5,4)(0,4)(0.5,4)(1,4)(1.5,4)(2,4)
\psdots(-0.5,5)(0,5)(0.5,5)(1,5)(1.5,5)(2,5)
\rput[c](3,4){$\cdots$}
\end{pspicture}
\caption{Spectral sequences $^{3}E, ^{4}E$.}
\end{figure}
\end{center}
Take $i:S\to M\Z_{p}$ a map representing the unit of $\pi_{0}{M\Z_{p}}$. This induces a map of spectra $K=K\wedge S\stackrel{1\wedge i}{\rightarrow} K\wedge M\Z_{p}$. This map induces a map of spectral sequences $j_{n}^{r,s}:\ ^{3}E_{n}^{r,s}\to\ ^{1}E_{n}^{r,s}$. Since each term of the spectral sequence $^{1}E$ is a $\Z_{p}$-module, by tensoring with $\Z_{p}$ we get a map of spectral sequences $\tilde{j}_{n}^{r,s}:\ ^{3}E_{n}^{r,s}\otimes \Z_{p}\to\ ^{1}E_{n}^{r,s}$. Notice that already on the $E_{2}$-level this map is an isomorphism because $H_{n}(BG)$ is finitely generated, and thus by  \cite[Corollary 56.4]{Munkres} we have a short exact sequence
\begin{equation}\label{eq15}
0\to H^{r}(BG,\Z)\otimes \Z_{p}\to H^{r}(BG,\Z_{p}) \to \text{Tor}_{1}^{\Z}(H^{r+1}(BG),\Z_{p})\to 0.
\end{equation}
Since $\Z_{p}$ is a flat $\Z$-module, from (\ref{eq15}) we see that 
\[
H^{r}(BG,\Z)\otimes \Z_{p}\approx H^{r}(BG,\Z_{p}),
\] 
and this isomorphism is precisely the $\tilde{j}$ map. Because the differentials with source $H^{2}(BG,K^{2}(*))$ in the spectral sequence $\{^{1}E_{n}^{r,s}\}$ are all trivial it follows that all the differentials with source $H^{2}(BG,K_{\Z_{p}}^{2}(*))$ are also trivial.
\qed

\begin{definition} Given a system of groups
\[
\{G_{n}\}=\cdots \to G_{n+1}\cdots \to G_{2}\to G_{1},
\] 
we say that $\{G_{n}\}$ satisfies the Mittag-Leffler condition if for every $i$ we can find a $j>i$ such that for every $k>j$ 
\[
\textsl{Im}(G_{k}\to G_{i})=\textsl{Im}(G_{j}\to G_{i}).
\]
\end{definition}

It is well known that if $\{G_{n}\}$ satisfies the Mittag-Leffler condition then 
\[
\limis_{k\to \infty}G_{k}=0.
\] 
On the other hand, if each $G_{k}$ is a countable group, then by \cite[Theorem 2]{McGibbon} we have that the system $\{G_{n}\}$ must satisfy the Mittag-Leffler condition.

Suppose now that $G$ is a compact Lie group. Then higher twistings of the Borel cohomology theory associated to $G$-equivariant K-theory and over a point are classified by 
\[
bsu_{\ot}^{1}(BG)=[BG,BBSU_{\otimes}].
\]
We are now able to show that for a compact Lie group this vanishes. We do this in the following theorem. 

\begin{theorem}
For any compact Lie group $G$, 
\[
bsu_{\otimes}^{1}(BG)=[BG,BBSU_{\otimes}]=0.
\]  
\end{theorem}
\Proof
For every $k \ge 0$ denote by $F_{k}$ the image of $\coprod_{0\le n\le k}{(G^{n}\times \Delta_{n})}$ in $BG$. The $F_{k}$'s form an increasing filtration of $BG$ and since $G$ is compact Lie each $F_{k}$ is of the homotopy type of a finite CW-complex. Let us denote 
\[
A_{k}=k^{4}(F_{k}) \text{\ \ and\ \ }  B_{k}=bsu_{\otimes}^{0}(F_{k}). 
\]
Using the filtration $\{F_{k}\}$ we get a short exact sequence 
\begin{equation}\label{eq16}
0\to \limis_{k\to \infty} A_{k}\to k^{5}(BG)\to \lim_{k\to \infty}k^{5}(F_{k})\to 0.
\end{equation}
By Theorem \ref{Bruner} we have that the middle term in (\ref{eq16}) vanishes and thus we see that $\limis_{k\to \infty} A_{k}=0$. By looking at the Atiyah-Hirzebruch spectral sequence, since $F_{k}$ is of the homotopy type of a finite CW-complex, we see that each $A_{k}$ and $B_{k}$ is finitely generated, in particular countable. Therefore the system $\{A_{k}\}$ satisfies the Mittag-Leffler condition.

On the other hand, by \cite[Corollary 1.4]{Adams} we have that after localization or completion at any prime $p$, the spectrum $bsu_{\otimes}$ is unique up to equivalence. In our context this means that $K\left\langle 4\right\rangle \wedge M\Z_{p}\simeq bsu_{\otimes}\wedge M\Z_{p}$ for every prime $p$. But we have that $bsu_{\otimes}\wedge M\Z_{p}\simeq \Sigma^{4} k\wedge M\Z_{p}$. Thus for each $k$ we have that 
\begin{equation}\label{eq17}
A_{k}\otimes \Z_{p}=k_{\Z_{p}}^{4}\stackrel{\simeq}{\rightarrow}(bsu_{\otimes}\wedge M\Z_{p})^{0}(F_{k})=B_{k}\otimes \Z_{p}.
\end{equation}
The outer equalities follow by \cite[Proposition 6.6]{Adamsbook}. Therefore we have a commutative diagram in which the vertical arrows are isomorphisms
\begin{equation}\label{chp7,1}
\begin{array}{ccccccc}
\rightarrow  & A_{n}\otimes \Z_{p} & \cdots  & \rightarrow  & A_{2}\otimes
\Z_{p} & \rightarrow  & A_{1}\otimes \Z_{p} \\ 
& \downarrow  &  &  & \downarrow  &  & \downarrow  \\ 
\rightarrow  & B_{n}\otimes \Z_{p} & \cdots  & \rightarrow  & B_{2}\otimes
\Z_{p} & \rightarrow  & B_{1}\otimes \Z_{p}.%
\end{array}%
\end{equation}%
Let $i>0$ be fixed. Since the system $\{A_{k}\}$ satisfies the Mittag-Leffler property we can find a $j>i$ such that for each $k>j$ 
\[\text{Im}(A_{k}\to A_{i})=\text{Im}(A_{j}\to A_{i}).
\]
The following is a short exact sequence:
\begin{equation}\label{eq19}
0\to \text{Ker}(A_{k}\to A_{i})\to A_{k}\to \text{Im}(A_{k}\to A_{i})\to 0.
\end{equation}
Since $\Z_{p}$ is a flat $\Z$-module we have that 
\begin{equation}\label{eq20}
0\to \text{Ker}(A_{k}\to A_{i})\otimes \Z_{p}\to A_{k}\otimes \Z_{p}\to \text{Im}(A_{k}\to A_{i})\otimes \Z_{p}\to 0
\end{equation}
is also exact. This shows that $\text{Im}(A_{k}\to A_{i})\otimes \Z_{p}=\text{Im}(A_{k}\otimes \Z_{p}\to A_{i}\otimes \Z_{p})$ and thus we see that for every $k>j$ and every prime $p$ we have 
\[
\text{Im}(A_{k}\otimes \Z_{p}\to A_{i}\otimes \Z_{p})=\text{Im}(A_{j}\otimes \Z_{p}\to A_{i}\otimes \Z_{p}).
\]
By the diagram (\ref{chp7,1}) we conclude that for every $p$ 
\begin{align*}
&\text{Im}(B_{k}\to B_{i})\otimes \Z_{p}=\text{Im}(B_{k}\otimes \Z_{p}\to B_{i}\otimes \Z_{p})=\\
&\text{Im}(B_{j}\otimes \Z_{p}\to B_{i}\otimes \Z_{p})=\text{Im}(B_{j}\to B_{i})\otimes \Z_{p}.
\end{align*}
Thus the groups $\text{Im}(B_{k}\to B_{i})$ and $\text{Im}(B_{j}\to B_{i})$ are two finitely generated groups that are equal after tensoring with $\Z_{p}$. By Lemma \ref{lemma 1} below we see that 
\[
\text{Im}(B_{k}\to B_{i})= \text{Im}(B_{j}\to B_{i}).
\]
We have proved that the system $\{B_{k}\}$ satisfies the Mittag-Leffler condition and thus
\begin{equation}\label{eq21}
\limis_{k\to \infty}B_{k}=\limis_{k\to \infty}bsu_{\otimes}^{0}(F_{k})=0.
\end{equation}
Using the filtration $\{F_{k}\}$ for the spectrum $bsu_{\otimes}$ we get a short exact sequence
\begin{equation}\label{eq22}
0\to \limis_{k\to \infty} B_{k}\to bsu_{\otimes}^{1}(BG)\to \lim_{k\to \infty}bsu_{\otimes}^{1}(F_{k})\to 0.
\end{equation}
Since the $\text{lim}^1$ part vanishes we get that
\[
bsu_{\otimes}^{1}(BG)=\lim_{k\to \infty}bsu_{\otimes}^{1}(F_{k}).
\]
We show now that the latter vanishes. To see this, note that for every prime $p$ we have a short exact sequence 
\begin{equation}\label{eq23}
0\to \limis_{k\to \infty}(bsu_{\otimes}\wedge M\Z_{p})^{0}(F_{k})\to (bsu_{\otimes}\wedge M\Z_{p})^{1}(BG)\to \lim_{k\to \infty}(bsu_{\otimes}\wedge M\Z_{p})^{1}(F_{k})\to 0.
\end{equation}
The term in the middle of (\ref{eq23}) vanishes and hence we see that 
\[
\lim_{k\to \infty}(bsu_{\otimes}\wedge M\Z_{p})^{1}(F_{k})=0.
\]
But by \cite[Proposition 6.6]{Adamsbook} we have that $(bsu_{\otimes}\wedge M\Z_{p})^{1}(F_{k})=bsu_{\otimes}^{1}(F_{k})\otimes \Z_{p}$. Thus for every prime $p$ the map 
\[
\lim_{k\to \infty}bsu_{\otimes}^{1}(F_{k})\otimes \Z_{p}=0.
\]
The proof finishes by using Lemma \ref{lemma 2} to see that 
\[
\lim_{k\to \infty}bsu_{\otimes}^{1}(F_{k})=0.
\]
\qed

\begin{lemma}\label{lemma 1}
Suppose that $A$ and $B$ are two finitely generated abelian groups with $A\subset B$ and that for every prime $p$, $A\otimes \Z_{p}=B\otimes \Z_{p}$. Then $A=B$.
\end{lemma}
\Proof
We have a short exact sequence
\[
0\to A\to B\to B/A\to 0.
\]
Since $\Z_{p}$ is a flat $\Z$-module we see that 
\[
0\to A\otimes \Z_{p}\to B\otimes \Z_{p}\to B/A\otimes \Z_{p}\to 0
\]
is also exact. As $A\otimes \Z_{p}= B\otimes \Z_{p}$ we see that $B/A\otimes \Z_{p}=0$. This is true for every $p$. This implies that $B/A=0$.
\qed

\begin{lemma}\label{lemma 2}
Let $\cdots \stackrel{f_{k}}{\rightarrow}G_{k}\stackrel{f_{k-1}}{\rightarrow} \cdots \stackrel{f_{2}}{\rightarrow} G_{2}\stackrel{f_{1}}{\rightarrow} G_{1}$ be a system of finitely generated abelian groups such that $\lim_{k\to \infty}G_{k}\otimes \Z_{p}=0$ for all primes $p$. Then $\lim_{k\to \infty}G_{k}=0$.
\end{lemma}
\Proof
Let $f:\prod_{i\ge 1}{G_{i}}\to \prod_{i\ge 1}{G_{i}}$ be defined by 
\[
f(x_{1},x_{2},...)=(x_{1}-f_{1}(x_{2}),x_{2}-f_{2}(x_{3}),...).
\]
We want to show that $f$ is injective, as $\lim_{k\to \infty}G_{k}=\text{Ker} (f)$. 
Suppose 
\[
x=(x_{1},x_{2},...)\in \text{Ker}(f).
\]
Then we have that $i_{p}(x)\in \text{Ker}(f_{p})=0$. Here $i_{p}:\prod_{i\ge 1}{G_{i}}\to \prod_{i\ge 1}{G_{i}\otimes \Z_{p}}$. Thus for each $i$ we have that $x_{i}\in \text{Ker}(G_{k}\to G_{k}\otimes \Z_{p})$ for each prime $p$. Since $G_{k}$ is finitely generated we have that 
\[
\bigcap_{p\ \ prime}\text{Ker}(G_{k}\to G_{k}\otimes \Z_{p})=0.
\] 
Thus $x=0$. 
\qed

\begin{corollary}
For a compact Lie grouop $G$ there are no higher twistings for the Borel cohomology theory associated to $G$-equivariant K-theory.
\end{corollary}
\medskip

\noindent{\bf{Remark:}} In general, the group $bsu_{\ot}^{1}(BG)$ does not vanish if $G$ does not satisfy the Atiyah-Segal Completion Theorem. To see this
let us consider an odd dimensional sphere ${S}^{2n+1}$ with $n\ge 2$. By the Kan-Thurston Theorem (see \cite{Kan-Thurston}) we know that there is a discrete group $G_{n}$ and a map $f:BG_{n}\to {S}^{2n+1}$ that is a homology equivalence. Since $f$ is a homology equivalence, it follows that $bsu_{\otimes}^{1}(BG_{n})=bsu_{\otimes}^{1}({S}^{2n+1})$. (This follows as we get isomorphism in the $E_{2}$-term and onward in the Atiyah-Hirzebruch spectral sequence.) Let us show now that  $bsu_{\otimes}^{1}({S}^{2n+1})\ne 0$. This will prove the proposition. Let $p$ be a prime number. We know that $bsu_{\otimes}\wedge M\Z_{p}\simeq bsu_{\oplus}\wedge M\Z_{p}\simeq \Sigma^{4} k\wedge M\Z_{p}$, and thus  
\[
bsu_{\otimes}^{1}({S}^{2n+1})\otimes \Z_{p}=(bsu_{\otimes}\wedge M\Z_{p})^{1}({S}^{2n+1})=k_{\Z_{p}}^{5}({S}^{2n+1})=k^{5}({S}^{2n+1})\otimes \Z_{p}.
\]
Here we used \cite[Proposition 6.6]{Adamsbook} as ${S}^{2n+1}$ is finite, and also the fact that $\Z_{p}$ is a flat $\Z$-module. Notice that both $bsu_{\otimes}^{1}({S}^{2n+1})$ and $k^{5}({S}^{2n+1})$ are finitely generated abelian groups. In general, if $A$ is a finitely generated abelian group, $A=0$ if and only if $A\otimes \Z_{p}=0$ for every prime number $p$. Thus, to show that $bsu_{\otimes}^{1}({S}^{2n+1})\ne 0$, we only need to show that $k^{5}(S^{2n+1})\ne 0$. To do so we use the Atiyah-Hirzebruch spectral sequence
\[
H^{r}(S^{2n+1},k^{s}(*))\Longrightarrow k^{r+s}(S^{2n+1}).
\]
\begin{center}
\begin{figure}[h]
\centering
\begin{pspicture}(-5,0)(5,6)
\psline{->}(-1,0)(-1,6)
\psline{->}(-1.5,3)(4.5,3)
\psline{->}(-1,3)(0,2.5)
\psdots(-1,1)
\psdots(-1,2)
\psdots(-1,3)
\psdots(2,1)
\psdots(2,2)
\psdots(2,3)
\rput[c](2,3.3){$_{2n+1}$}
\rput[c](4.5,2.8){$r$}\rput[c](-1.3,6){$s$}\rput[c](4.5,6){$E_{2}^{r,s}$}
\rput[c](0.5,0.5){$\vdots$}
\end{pspicture}
\caption{Atiyah-Hirzebruch spectral sequences for $ k^{r+s}(S^{2n+1})$.}
\end{figure}
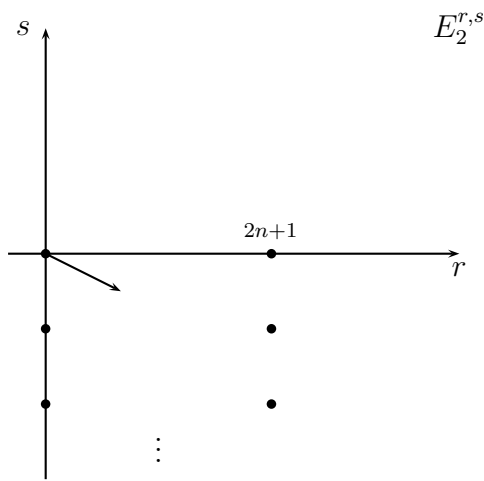
\end{center}
We claim that this spectral sequence collapses on the $E_{2}$-term. To see this, we only need to note that the corresponding spectral sequence collapses in the case of $K$. Since we have a map of spectra $k\to K$ inducing an isomorphism on $\pi_{n}$ for $n\ge0$, the spectral sequence in the case of $k$ also collapses. Since  $n\ge 2$, we see $k^{5}(S^{2n+1})=\Z\ne0$.
\medskip

\noindent{\bf{Remark:}} If $G$ is a compact Lie group and if we consider twistings of the Borel cohomology theory associated with $G$-equivariant K-theory we encounter higher twistings if we work with spaces more general than a point. For example in the trivial case where $G=\{e\}$ is the trivial group, then for $X=S^{2n+1}$ an odd sphere with $n\ge 2$ we have higher twistings these are classified by the group
\[
bsu_{\ot}(S^{2n+1})=\Z\ne 0.
\]


\begin{thebibliography}{99}

\bibitem{Adamsbook} J. F. Adams. Stable homotopy and generalized homology. Univ. of Chicago Press, 1974.

\bibitem{Adams} J. F. Adams and S. B. Priddy. Uniqueness of BSO. Math. Proc. Camb. Soc. (1976), 80-475.

\bibitem{Atiyahcompletion} M. F. Atiyah, and G. Segal. Equivariant K-theory and completion, J. Differential Geom. 3 (1969), 1–18.

\bibitem{Atiyah} M. F. Atiyah and G. Segal. Twisted K-theory. Ukr. Mat. Visn.  1  (2004),  no. 3, 287--330;  translation in  Ukr. Math. Bull.  1  (2004),  no. 3, 291--334
55N15 (19Kxx 46L80 55N91)

\bibitem{Atiyahtwisted} M. F. Atiyah, and G. Segal. Twisted K-theory and cohomology, Arxiv e-print.

\bibitem{Kan-Thurston} D. M. Kan. W. P. Thurston. Every connected space has the homology of a $K(\pi,1)$.  Topology  15  (1976), no. 3, 253-258. 

\bibitem{Madsen} I. Madsen, V. Snaith and J. Tornehave. Infinite loop maps in geometric topology.  Math. Proc. Cambridge Philos. Soc.  81  (1977), no. 3, 399--430. 

\bibitem{Mayring} J. P. May. $E_{\infty}$ ring spaces and $E_{\infty}$ ring spectra. With contributions by Frank Quinn, Nigel Ray, and Jorgen Tornehave. Lecture Notes in Mathematics, Vol. 577. Springer-Verlag, Berlin-New York, 1977.

\bibitem{MayEinfty} J. P. May. What precisely are {$E\sb{\infty }$} ring spaces and {$E\sb{\infty }$} ring  spectra?, Preprint 2008.
          
\bibitem{McGibbon} C. A. McGibbon and J. M. Moller. On spaces with the same n-type for all n. Topology 31 (1992).

\bibitem{Munkres} J.R. Munkres. Elements of Algebraic Topology. Perseus Publishing. Cambridge. 1984.

\bibitem{Segal} G. Segal, Categories and cohomology theories, Topology, 13, (1974), 293-312


\end{thebibliography}
\end{document}